\documentclass[fleqn]{mat01}
\usepackage{times,mathtimy,amssymb,latexsym}
\begin{document}

\def\d{\mbox{\rm d}}
\def\e{\mbox{\rm e}}

\newtheorem{theore}{Theorem}
\renewcommand\thetheore{\arabic{section}.\arabic{theore}}
\newtheorem{theor}[theore]{\bf Theorem}
\newtheorem{coro}[theore]{\rm COROLLARY}
\renewcommand\theequation{\arabic{section}.\arabic{equation}}

\setcounter{page}{371}
\firstpage{371}

\title{On the series $\pmb{\sum_{k = 1}^\infty \binom{3k}{k}^{ - 1}k^{ - n}\,x^k}$}

\markboth{Necdet Batir}{On the series $\sum_{k = 1}^\infty \binom{3k}{k}^{ - 1}k^{ - n}\,x^k$}

\author{NECDET BATIR}

\address{Department of Mathematics, Faculty of Arts and Sciences,
Y\"{u}z\"{u}nc\"{u} Yil University, 65080 Van, Turkey\\
\noindent E-mail: necdet{\_}batir@hotmail.com}

\volume{115}

\mon{November}

\parts{4}

\pubyear{2005}

\Date{MS received 6 October 2004; revised 20 August 2005}

\begin{abstract}
In this paper we investigate the series $\sum_{k =
1}^\infty \binom{3k}{k}^{ - 1}k^{ - n}x^k$. Obtaining some integral
representations of them, we evaluated the sum of them explicitly for
$n = 0, 1, 2$.
\end{abstract}


\maketitle

\section{Introduction}

After Ap\'{e}ry \cite{2} proved the irrationality of $\zeta (2)$ and $\zeta (3)$,
where $\zeta $ is the Riemann-zeta function defined by
\begin{equation*}
\zeta (s) = \sum\limits_{k = 1}^\infty {\frac{1}{k^s}}, \quad \hbox{Re}\,s\, > 1,
\end{equation*}
by employing the series representations
\begin{equation*}
\zeta (2) = 3\sum\limits_{n = 1}^\infty \frac{1}{n^2 \binom{2n}{n}}
\quad \hbox{and} \quad \zeta (3) = \frac{5}{2}\sum\limits_{n =
1}^\infty \frac{( - 1)^{n - 1}}{n^3 \binom{2n}{n}},
\end{equation*}
many authors have considered the series involving inverse binomial coefficients
and they obtained many interesting results. We have a similar series
representation for $\zeta (4)$:
\begin{equation*}
\zeta (4) = \frac{36}{17}\sum\limits_{n = 1}^\infty \frac{1}{n^4
\binom{2n}{n}},
\end{equation*}
see \cite{8}. Some other related interesting results involving binomial
coefficients can be found in Chapter~9 of \cite{4}, \cite{3},
\cite{5,6,7,8,9} and \cite{11,12}.

Motivated by such results we shall consider here the following family of
sums:
\begin{equation*}
S(n,m;x) = \sum\limits_{k = 1}^\infty \frac{x^k}{k^n \binom{3mk}{mk}}.
\end{equation*}
A good way to approach these series is to try and find their integral
representations. In this way we can evaluate many of them explicitly.

In this paper, we will use, as usual, the following definitions and
identities for the Euler's gamma function $\Gamma$, beta function $\beta$,
polylogarithms $Li_n (z)$ and generalized hypergeometric series ${ }_pF_q
(a_1, a_2, \dots, a_p; b_1, b_2, \dots, b_q\hbox{\rm :}\ x)$,
\begin{equation}
\beta (s, t) = \int^{1}_{0} u^{s - 1} (1 - u)^{t - 1} \d u =
\frac{\Gamma (s)\Gamma (t)}{\Gamma (s + t)} \quad \hbox{for} \ \ s > 0, t > 0,
\end{equation}
(see Theorem~7.69 of \cite{4}),
\begin{align}
Li_n (z) &:= \sum\limits_{{\kern 1pt} k = 1}^{{\kern 1pt} \infty }
{\frac{z^k}{k^n}} = \frac{( - 1)^{n - 1}}{(n -
1)!}\int_{0}^{1} {\frac{z \log^{n -
1}\phi \d \phi}{1 - z \phi}} \quad \hbox{for} \ \ |z| \leq 1,\\[.2pc]
Li_n (z^m) &= m^{n - 1} \sum\limits_{k = 1}^m {Li_n (\omega^kz)},
\end{align}
where $m$ is a positive integer and $\omega = \hbox{e}^{2\pi i/m}$. The $m$th primitive
root of unity is called factorization formula for polylogarithm series and
\begin{align*}
{ }_pF_q (a_1, a_2, \dots, a_p; b_1, b_2, \dots, b_q\hbox{\rm :}\ x) =
\sum\limits_{k = 1}^\infty {\frac{(a_1 )_k (a_2 )_k \dots (a_p )_k x^k}{(b_1
)_k (b_2 )_k \dots (b_q )_k }},
\end{align*}
where
\begin{equation*}
(a)_k = \frac{\Gamma (a + k)}{\Gamma (a)}.
\end{equation*}
For further properties of polylogarithms and hypergeometric series and
related functions, see \cite{10} and Chapter~2 of \cite{1}, respectively. Almost all
results given here were obtained using identities (1.1) and (1.2)
extensively.

\section{Main results}

The main results of this paper are the following two theorems.

\begin{theor}[\!]
For $|x| \leq 27/4$ and $n = 2, 3, \dots$ we have
\setcounter{equation}{0}
\begin{align}
\sum\limits_{k = 1}^\infty \frac{x^k}{k^n \binom{3k}{k}} &= \frac{(- 1)^{n - 1}}{(n - 2)!}\int_{0}^{\alpha (x)} {u\log ^{n
- 2}\left[ {\frac{1}{x}\frac{(1 - \e^u)^3}{\e^{2u}}} \right]\,} \d
u\nonumber\\[.2pc]
&\quad\, +\!\frac{4( - 1)^{n - 2}}{3(n - 2)!}\int_{0}^{\beta
(x)}\!\!{v\log ^{n - 2}\!\left[\!{\frac{[1\!+\!2 \cos[(2v\!+\!2\pi)/3]]^3}{2x[1\!+\!\cos [(2v\!+\!2\pi)/3]]}}\!\right]\d v},
\end{align}
where
\begin{equation*}
\alpha (x) = \log \left[ {\frac{\phi ^3 + 1}{[\phi + 1]^3}} \right], \quad
\beta (x) = 3 \arctan \left[ {\frac{\sqrt 3 }{2\phi - 1}} \right]
\end{equation*}
and
\begin{equation*}
\phi (x):= \left[ {\frac{27 - 2x + 3 [81 - 12x]^{1/2}}{2x}} \right]^{1/3}.
\end{equation*}
\end{theor}

\begin{proof}
We start with identity (1.1).
\begin{align*}
\sum\limits_{k = 1}^\infty \frac{x^k}{k^n \binom{3k}{k}} &= \sum\limits_{k = 1}^\infty
{\frac{x^k(k!)(2k)!}{k^n(3k)!}}\\[.2pc]
&= \sum\limits_{k = 1}^\infty {\frac{x^k\Gamma (k + 1)\Gamma (2k +
1)}{k^n\Gamma (3k + 1)}}\\[.2pc]
&= \sum\limits_{k = 1}^\infty {\frac{x^k\Gamma (k)\Gamma (2k + 1)}{k^{n -
1}\Gamma (3k + 1)}}\\[.2pc]
&= \sum\limits_{k = 1}^\infty {\frac{x^k}{k^{n - 1}}} \beta (k, 2k +
1)\\[.2pc]
&= \sum\limits_{k = 1}^\infty {\frac{x^k}{k^{n - 1}}} \int_{0}^{1}
{t^{k - 1}(1 - t)^{2k} \d t}.
\end{align*}
Inverting the order of summation and integration, we get
\begin{align}
\sum\limits_{k = 1}^\infty \frac{x^k}{k^n \binom{3k}{k}} &=
\int_{0}^{1} {\sum\limits_{k = 1}^\infty {\frac{[xt(1 -
t)^2]^k}{k^{n - 1}}} \frac{\d t}{t}}\nonumber\\[.2pc]
&= \int_{0}^{1} {\frac{Li_{n - 1} [xt(1 - t)^2]}{t}} \d t\\[.2pc]
&= \frac{( - 1)^{n - 2}}{(n - 2)!}\int_{0}^{1}\!{\left[\!
{\int_{0}^{1} {\frac{xt(1 - t)^2\log ^{n - 2}z}{1 - xt(1 -
t)^2z}\d z}} \!\right]}\!\frac{\d t}{t},
\end{align}
where in the last step we employ identity (1.2). Inverting the order of
integration here and leaving the justification of it at the end of the proof,
we obtain
\begin{align*}
\sum\limits_{k = 1}^\infty \frac{x^k}{k^n \binom{3k}{k}} = \frac{( - 1)^{n - 1}}{(n - 2)!}\int_{0}^{1} {\frac{\log ^{n -
1}z}{z}\left[\!{\int_{0}^{1} {\frac{t^2 - 2t + 1}{t^3 - 2t^2 + t - [xz]^{
- 1}}\d t} }\!\right]} \d z.
\end{align*}
Making the change of variable $t = u + 2/3$ here, we find after some
manipulations that
\begin{align*}
\hskip -4pc \sum\limits_{k = 1}^\infty \frac{x^k}{k^n \binom{3k}{k}} = \frac{( - 1)^{n - 2}}{3(n - 2)!}\int_{0}^{1} {\frac{\log ^{n -
2}z}{z}\left[\!{\int_{ - 2/3}^{1/3} {\frac{2u - 2/3}{u^3 - u/3
+ 2/27 - [xz]^{ - 1}}\d u} }\!\right]} \d z.
\end{align*}
Now making the change of variable
\begin{equation*}
z = \frac{27}{x}\frac{y^3}{(y^3 + 1)^2}
\end{equation*}
in the first integral, we obtain
\begin{align*}
\sum\limits_{k = 1}^\infty \frac{x^k}{k^n \binom{3k}{k}} &= \frac{( - 1)^{n - 1}}{(n - 2)!}\int_{0}^{\phi (x)} {\log ^{n -
2}\left[ {\frac{27}{x}\frac{y^3}{(y^3 + 1)^2}} \right]}\\[.2pc]
&\quad\,\times \left[ {\int_{- 2/3}^{1/3} {\frac{2u - 2/3}{u^3 - u/3 -
(1 + y^6)/(27y^3)}\d u}} \right]\frac{1}{y}\frac{y^3 - 1}{y^3 + 1}\d y.
\end{align*}
Here
\begin{equation}
\phi (x) = \left[ {\frac{27 - 2x + 3 [81 - 12x]^{1/2}}{2x}} \right]^{1/3}.
\end{equation}
If we make the change of variable $t = 3y/(y^2 + 1)$ in this integral, we
find that
\begin{align}
\hskip -4pc \sum\limits_{k = 1}^\infty \frac{x^k}{k^n \binom{3k}{k}} &= \frac{3( - 1)^{n - 1}}{(n - 2)!}\int_{0}^{\lambda (x)} {\log ^{n -
2} \left[ {\frac{27}{x}} \frac{t^3}{(2t + 3)(t - 3)^2} \right]}\nonumber\\[.2pc]
\hskip -4pc &\quad\, \times \left[ {\int_{- 2/3}^{1/3} {\frac{2u - 2/3}{u^3 - u/3
- (3 - t^2)/3t^3}\,\d u} } \right]\,\frac{t + 3}{t(t - 3)(2t + 3)}\d t,
\end{align}
where
\begin{equation*}
\lambda (x) = \frac{3\phi(x)}{\phi (x)^2 + 1}
\end{equation*}
with $\phi (x)$ defined by (2.4). First, we compute the inner integral.
By Cardano's method, the roots of the cubic equation $u^3 - u/3 - (3 -
t^2)/3t^3 = 0$ are
\begin{equation*}
\hskip -4pc \alpha = 1/t, \quad \beta = [ - 3 - i\sqrt{27 - 12t^2}]/6t \quad \hbox{and} \quad \gamma =
 - 3 + i\sqrt{27 - 12t^2}]/6t.
\end{equation*}
Thus, we can factorize the integrand in the inner integral as
\begin{align*}
\frac{2u - 2/3}{u^3 - u/3 - (3 - t^2)/3t^3} &= \frac{2t}{t +
3}\frac{1}{u - 1/t} - \frac{t}{t + 3}\frac{2u + 1/t}{u^2 + u/t + 1/t^2 -
1/3}\\[.2pc]
&\quad\, - \frac{2t + 3}{t + 3}\frac{1}{u^2 + u/t + 1/t^2 - 1/3}.
\end{align*}
Integrating both sides of this equation from $-2/3$ to 1/3 and then
simplifying it, we find that
\begin{align*}
&\int_{- 2/3}^{1/3} {\frac{2u - 2/3}{u^3 - u/3 - (3 - t^2)/3t^3}\,\d
u}\\
&\quad\, = \frac{3t}{t + 3}\log \left[ {\frac{3 - t}{3 + 2t}} \right] - \frac{2t\!+\!3}{t\!+\!3}\frac{6t}{[27\!-
\!12t^2]^{1/2}}\arctan\!\left[\!\!{\frac{3t}{5t\!-\!6}\sqrt {\frac{9\!-
\!6t}{3\!+\!2t}}}\!\right]\!.
\end{align*}
Replacing this in (2.5), we obtain after some simplification
\begin{equation}
\sum\limits_{k = 1}^\infty \frac{x^k}{k^n \binom{3k}{k}} = S_1 + S_2,
\end{equation}
where
\begin{equation}
\hskip -4pc S_1 = \frac{9( - 1)^{n - 1}}{(n - 2)!}\int_{0}^{\lambda (x)} {\log ^{n
- 2}\left[ {\frac{27}{x}\frac{t^3}{(2t + 3)(t - 3)^2}} \right]} \log
\left[ {\frac{3 - t}{2t + 3}} \right]\,\frac{\d t}{(t - 3)(2t + 3)}
\end{equation}
and
\begin{align}
S_2 &= \frac{18( - 1)^{n - 2}}{(n - 2)!}\int_{0}^{\lambda (x)}
\log^{n - 2}\left[ \frac{27}{x}\frac{t^3}{(2t + 3)(t - 3)^2} \right]\nonumber\\[.2pc]
&\quad\, \times \arctan \left[ \frac{3t}{5t - 6}\sqrt \frac{9 - 6t}{3 + 2t}
\right]\frac{\d t}{(t - 3)[27 - 12t^2]^{1/2}}.
\end{align}
Now we simplify these two integrals. If we make in (2.7) the change of variable
\begin{equation*}
u = \log \left[ {\frac{3 - t}{2t + 3}} \right],
\end{equation*}
we find that
\begin{equation}
S_1 = \frac{( - 1)^{n - 1}}{(n - 2)!}\int_{0}^{\alpha (x)} {u\log ^{n
- 2}\left[ {\frac{1}{x}\frac{(1 - \e^u)^3}{\e^{2u}}} \right]\,} \d u,
\end{equation}
where
\begin{equation*}
\alpha (x) = \log \left[ {\frac{3 - \lambda (x)}{3 + 2\lambda (x)}}
\right].
\end{equation*}
In (2.8), making the change of variable
\begin{equation*}
y = \sqrt {\frac{9 - 6t}{3 + 2t}},
\end{equation*}
we arrive at the following:
\begin{equation}
S_2 = \frac{4( - 1)^{n - 2}}{(n - 2)!}\int_{\sqrt 3 }^{\lambda _1
(x)} {\log ^{n - 2}\left[ {\frac{(3 - y^2)^3}{4x(y^2 + 1)^2}}
\right]\,\frac{(3\arctan y - \pi )}{y^2 + 1}\,\d y}
\end{equation}
since
\begin{equation*}
\arctan \left[ {\frac{y^3 - 3y}{3y^2 - 1}} \right] = 3 \arctan y
- \pi, \quad \hbox{for} \ \ y > 0
\end{equation*}
where
\begin{equation*}
\lambda _1 (x) = \sqrt {\frac{9 - 6\lambda (x)}{3 + 2\lambda (x)}}.
\end{equation*}
We need to induce one more change of variable to bring (2.10) in a simple
form. Setting $v = 3 \arctan y - \pi $ here, we get
\begin{equation}
S_2 = \frac{4( - 1)^{n - 2}}{3(n - 2)!}\int_{0}^{\beta (x)}
{v\log ^{n - 2}\left[ {\frac{[1 + 2\cos [(2v + 2\pi
)/3]]^3}{2x[1 + \cos [(2v + 2\pi)/3]]}} \right]\,\d v},
\end{equation}
\pagebreak

\noindent where
\begin{align*}
\beta (x) &= 3\arctan \sqrt {\frac{9 - 6\lambda (x)}{3 + 2\lambda (x)}} - \pi
= 3 \arctan \frac{\sqrt 3 |{\phi (x) - 1}|}{\phi (x) + 1} - \pi\\[.2pc]
&= 3\arctan \left[ {\frac{\sqrt 3 }{1 - 2\phi (x)}} \right].
\end{align*}
Substituting the values of $S_1$ and $S_2$ from (2.9) and (2.11) in
(2.6), we get
\begin{align*}
\hskip -4pc \sum\limits_{k = 1}^\infty \frac{x^k}{k^n \binom{3k}{k}} &= \frac{( - 1)^{n - 1}}{(n - 2)!}\int_{0}^{\alpha (x)} {u\log ^{n -
2}\left[ {\frac{1}{x}\frac{(1 - \e^u)^3}{\e^{2u}}} \right]\,} \d u\\[.2pc]
&\quad\, + \frac{4( - 1)^{n - 2}}{3(n - 2)!}\int_{0}^{\beta (x)}
{v\log ^{n - 2}\left[ {\frac{[1 + 2\cos [(2v + 2\pi
)/3]]^3}{2x[1 + \cos [(2v + 2\pi )/3]]}} \right]\,\d v},
\end{align*}
where
\begin{equation*}
\alpha (x) = \log \left[ {\frac{\phi (x)^3 + 1}{(\phi (x) + 1)^3}}
\right] \quad \hbox{and} \quad \beta (x) = 3\arctan \left[ {\frac{\sqrt 3}{1 - 2\phi (x)}}
\right].
\end{equation*}
To complete the proof of Theorem~2.1 we need to justify the inversion made
in (2.3). In the inner integral in (2.3), we induce the change of
variable $z = 1/u$ to get
\begin{equation*}
\hskip -4pc \int_{0}^{1} {\;\left[ {\int_{0}^{1} {\frac{xt(1 - t)^2\log ^{n -
2}z}{1 - xt(1 - t)^2z}\,\d z} } \right]} \,\d t = \int_{0}^{1}
{\;\left[ {\int_{1}^{\infty } {\frac{x(1 - t)^2\log ^{n
- 2}u}{u^2 - xt(1 - t)^2u}\,\d u} } \right]} \,\d t.
\end{equation*}
Since for every $0 \leq t \leq 1, - 27/4 \leq x \leq 27/4$
and $u \geq 1$,
\begin{equation*}
\frac{(1 - t)^2\log ^{n - 2}u}{u^2 - xt(1 - t)^2u} \le \frac{\log ^{n -
2}u}{u^2 - u}
\end{equation*}
and the improper integral
\begin{equation*}
\int_{1}^{\infty } {\frac{\log ^{n - 2}u}{u^2 - u}\,\d u}
\end{equation*}
is convergent, and
\begin{equation*}
\int_{1}^{\infty } {\frac{x(1 - t)^2\log ^{n - 2}u}{u^2 - xt(1 -
t)^2u}\,\d u}
\end{equation*}
is uniformly convergent. This justifies the inversion of the order of the
integrals in (2.3) and hence the proof of Theorem~2.1 is complete. \hfill $\Box$
\end{proof}

The next theorem gives a generalization of Theorem~2.1.

\begin{theor}[\!]
For $m = 1,2,3,\dots,$ $n = 1,2,3,\dots$ and $|x| \leq (27 /
4)^m,$ we have
\begin{equation}
S(n,m\hbox{\rm :}\ x) = m^{n - 1}\sum_{j = 1}^m S(n,1\hbox{\rm :}\
\omega^{j}x^{1 / m}){\rm ,}
\end{equation}
where $\omega = \e^{2\pi i / m}$ is a primitive root of unity.
\end{theor}

\begin{proof}
\begin{align*}
S(n,m\hbox{\rm :}\ x) &= m\sum_{k = 1}^\infty \frac{x^k\Gamma
(mk)\Gamma (2mk + 1)}{k^{n - 1}\Gamma (3mk + 1)}\\[.4pc]
&= m\sum_{k = 1}^\infty \frac{x^k}{k^{n - 1}}\beta
(mk,2mk + 1)\\[.4pc]
&= m\sum_{k = 1}^\infty \frac{x^k}{k^{n - 1}}\int_{0}^{1}
{t^{mk - 1}(1 - t)^{2mk}} \d t.
\end{align*}
Inverting the order of summation and integration, we find that
\begin{align*}
S(n,m\hbox{\rm :}\ x) &= m\int_{0}^{1} \sum_{k = 1}^\infty
\frac{[(t(1 - t)^2x^{1 / m})^m]^k}{k^{n - 1}} \frac{\d
t}{t}\\[.4pc]
&= m\int_{0}^{1} {\frac{Li_{n - 1} [(t(1 - t)^2x^{1 /
m})^m]}{t}} \d t\\[.4pc]
&= m^{n - 1}\sum_{j = 1}^m {\int_{0}^{1} {\frac{Li_{n -
1} (\omega ^jt(1 - t)^2x^{1 / m})}{t}}} \d t\\[.4pc]
&= m^{n - 1}\sum_{j = 1}^m {1\int_0^1 {\sum_{k = 1}^\infty
{\frac{[\omega ^jt(1 - t)^2x^{1 / m}]^k}{k^{n - 1}}} \frac{\d
t}{t}}}.
\end{align*}
Inverting the order of summation and integration, we get
\begin{align*}
S(n,m\hbox{\rm :}\ x) &= m^{n - 1}\sum_{j = 1}^m {\sum_{k =
1}^\infty {\frac{[\omega ^{jk}x^{1 / m}]^k}{k^{n - 1}}} \int_0^1
{t^{k - 1}(1 - t)^{2k}\d t}
}\\[.4pc]
&= m^{n - 1}\sum_{j = 1}^m \sum_{k = 1}^\infty \frac{[\omega
^jx^{1 / m}]^k}{k^n\binom{3k}{k}}\\[.4pc]
&= m^{n - 1}\sum_{j = 1}^m {S(n,1\hbox{\rm :}\ \omega ^jx^{1 /
m})},
\end{align*}
completing the proof of Theorem~2.2. \hfill $\Box$
\end{proof}

\begin{coro}$\left.\right.$\vspace{.5pc}

\noindent For $m = 1,2,3,\dots$ and $|x| \leq (27 / 4)^m$ we have
\begin{align}
\sum_{k = 1}^\infty \frac{x^k}{k^2\binom{3mk}{mk}} &= m\sum_{k =
1}^m \left\lbrace 6\arctan^{2}\left[\frac{\sqrt 3 }{2\phi (\omega
^kx^{1 / m}) -
1}\right]\right. \nonumber\\[.3pc]
&\quad\, \left. - \frac{1}{2}\log ^2\left[\frac{1 + [\phi (\omega
^kx^{1 / m})]^3}{[1 + \phi (\omega ^kx^{1 / m})]^3} \right]
\right\rbrace,
\end{align}
where $\omega = \e^{2\pi i / m}$ is a primitive root of unity.
\end{coro}

\begin{proof}
Setting $n = 2$ in (2.12) we get the desired result. \hfill $\Box$
\end{proof}

\section{Applications}

Putting some particuler values for $n$ and $x$ in Theorems~2.1 and
2.2, we can make many explicit evaluations.

Let $|x|\leq 27/4$. If we set $n = 2$ in (2.1) we find by the help
of Gauss multiplication formula for Euler's gamma function:
\setcounter{equation}{0}
\begin{align}
S(2,1\hbox{\rm :}\ x) &= \sum_{k = 1}^\infty \frac{x^k}{k^2\binom
{3k}{k}} = \frac{x}{3}\,{ }_4 F_3 \left(1,1, 1,\frac{3}{2};
\frac{4}{3},\frac{5}{3}, 2;
\frac{4x}{27}\right)\nonumber\\[.4pc]
&= 6\arctan^2\left[{\frac{\sqrt 3}{2\phi - 1}} \right] -
\frac{1}{2}\log^2\left[{\frac{\phi ^3 + 1}{(\phi + 1)^3}} \right].
\end{align}
Differentiating (3.1) with respect to $x$ and then multiplying by
$x$ we get for $|x| < 27/4$:
\begin{align}
S(1,1\hbox{\rm :}\ x) &= \sum_{k = 1}^\infty \frac{x^k}{k\binom
{3k} {k}} = \frac{x}{3}\,{ }_3 F_2 \left(1,1,\frac{3} {2};
\frac{4}{3}, \frac{5}{3};
\frac{4x}{27}\right)\nonumber\\[.4pc]
&= \frac{1}{\sqrt {27 - 4x} }\left\lbrace \arctan \left[
\frac{\sqrt 3}{2\phi - 1} \right]\frac{18\phi}{1 - \phi + \phi ^2}\right.\nonumber\\[.4pc]
&\quad\, \left. - \log \left[\frac{\phi ^3 + 1}{(\phi + 1)^3}
\right]\frac{3\sqrt 3 \phi (1 - \phi )}{1 + \phi ^3}
\right\rbrace.
\end{align}
Differentiating both sides of (3.2) with respect to $x$ and then
multiplying by $x$ we get for $|x| < 27/4$:
\begin{align}
\sum_{k = 1}^\infty \frac{x^k}{\binom {3k} {k}} &= { }_3 F_2
\left( {1, \frac{3}{2}, 2; \frac{4}{3}, \frac{5}{3}:
\frac{4x}{27}} \right)\nonumber\\[.4pc]
&= \left[ {\frac{36\phi x}{(27 - 4x)^{3 / 2}(1 - \phi + \phi ^2)}
- \frac{18\sqrt 3 (1 - \phi ^2)\phi }{(1 - \phi + \phi
^2)^2(27 - 4x)}} \right]\nonumber\\[.4pc]
&\quad\ \times \arctan \left[ {\frac{\sqrt 3 }{2\phi
- 1}}\right]\nonumber\\[.4pc]
&\quad  + \left[ {\frac{9\phi (1 - 2\phi - 2\phi ^3 + \phi
^4)}{(1 + \phi ^3)^2(27 - 4x)} - \frac{6\sqrt 3 (1 - \phi
)\phi x}{(27 - 4x)^{3 / 2}(1 + \phi ^3)}} \right]\nonumber\\[.3pc]
&\quad\,\, \times \log \left[ {\frac{1 + \phi ^3}{(1 + \phi )^3}}
\right] + \frac{108\phi ^3}{(27 - 4x)(1 + \phi ^3)^2},
\end{align}
where
\begin{equation*}
\phi = \phi (x) = \left[ {\frac{27 - 2x + 3 [81 - 12x]^{1 /
2}}{2x}} \right]^{1 / 3},
\end{equation*}
as defined by (2.4).

Putting $x = 27/4$ and $n = 2$ in (2.1) yields
\begin{equation}
\sum_{k = 1}^\infty \frac{(27 / 4)^k}{k^2\binom {3k} {k}} =
\frac{2\pi ^2}{3} - 2\log^22.
\end{equation}
Let $n = 2$ and $x = 6$ in (2.1). Then
\begin{equation}
\sum_{k = 1}^\infty \frac{6^k}{k^2\binom {3k} {k} } =
6\arctan^2\left[ {\frac{\sqrt 3 }{2^{4 / 3} - 1}} \right]
 - \frac{1}{2}\log ^2(2^{1 / 3} - 1).
\end{equation}
Set $x = 1/2$ in (3.1) to get
\begin{equation}
\sum_{k = 1}^\infty {\frac{1}{k^2\binom {3k} {k}2^k} =
\frac{1}{24}\pi ^2 - \frac{1}{2}\log ^22}.
\end{equation}
Set $x = 1$ in (3.1) to get
\begin{align}
\sum_{k = 1}^\infty \frac{1}{k^2\binom{3k} {k}} &= 6\arctan
^2\left[ {\frac{\sqrt 3 }{1 - [100 + 12\sqrt {69} ]^{1 / 3}}}
\right]\nonumber\\[.3pc]
&\quad\, - \frac{1}{2}\log ^2\left[ {\frac{12(9 + \sqrt
{69})}{[2 + (100 + 12\sqrt {69})^{1 / 3}]^3}} \right].
\end{align}
Set $x = -1/4$ in (3.1) to get
\begin{equation}
\sum_{k = 1}^\infty \frac{( - 1 / 4)^k}{k^2\binom {3k} {k}} =
6\,\hbox{arc}\cot ^2( {2\sqrt 3 + \sqrt 7 }) - \frac{1}{2}\log^22.
\end{equation}
Put $x = 1/2$ in (3.2) to get
\begin{equation}
\sum_{k = 1}^\infty \frac{1}{k\binom {3k} {k} 2^k} =
\frac{1}{10}\pi - \frac{1}{5}\log 2 .
\end{equation}
Set $x = 6$ in (3.2) to get
\begin{align}
\sum_{k = 1}^\infty \frac{6^k}{k\binom {3k} {k}} &= \sqrt 3 2^{4 /
3}(1 + 2^{1 / 3})\arctan \left[ {\frac{\sqrt 3 }{2^{4 / 3} - 1}}
\right]\nonumber\\[.3pc]
&\quad\,- 2^{1 / 3}(1 - 2^{1 / 3})\log (2^{1 / 3} - 1).
\end{align}
Let $x = 1/2$ and $n = 0$ in (3.3) to get
\begin{equation}
\sum_{k = 1}^\infty \frac{1}{\binom {3k} {k}2^k} = \frac{2}{25} -
\frac{6}{125}\log 2 + \frac{11}{250} \pi.
\end{equation}
Here, observe that $\phi (1 / 2) = 2 + \sqrt 3$. Let $x = -1/4$ in
(3.3) to get
\begin{equation}
\sum_{k = 1}^\infty {\frac{( - 1)^k}{\binom {3k} {k}4^k} = -
\frac{1}{28} - \frac{3}{32}\log 2 + \frac{39}{112\sqrt 7 }}
\hbox{arc}\cot (2\sqrt 3 + \sqrt 7).
\end{equation}
Note that $\phi ( - 1 / 4)= - (5 + \sqrt {21} ) / 2$. Set $x = 6$
in (3.3) to get
\begin{align}
\sum_{k = 1}^\infty \frac{6^k}{\binom{3k}{k}} &= 2\left( {240 +
96.2^{1 / 3} + 75.2^{2 / 3}} \right)^{1 / 2} \arctan \left(
{\frac{\sqrt 3 }{2^{4 / 3} - 1}} \right)\nonumber\\[.3pc]
&\quad\,+ 2^{1/3} ( {4.2^{1 / 3} - 5})\log (2^{1 /
3} - 1) + 8.
\end{align}
Setting $x = 1$ in (3.3) we obtain
\begin{align}
\sum_{k = 1}^\infty \frac{1}{\binom {3k}{k}} &= \left[
{\frac{36\sqrt {23} \tau }{529(1 - \tau + \tau ^2)} -
\frac{18\sqrt 3 (1 - \tau ^2)\tau }{23(1 - \tau + \tau ^2)^2}}
\right]\arctan \left[ {\frac{\sqrt 3 }{2\tau - 1}} \right]\nonumber\\[.4pc]
&\quad + \left[ {\frac{9\tau (1 - 2\tau - 2\tau ^3 + \tau
^4)}{23(1 + \tau ^3)^2} - \frac{6\sqrt {69} (1 - \tau )\tau
}{529(1 + \tau ^3)}} \right]\nonumber\\[.3pc]
&\qquad\, \times \log \left[ {\frac{1 + \tau ^3}{(1 + \tau )^3}}
\right] + \frac{108\tau ^3}{23(1 + \tau ^3)^2},
\end{align}
where
\begin{equation*}
\tau = \left(\frac{25 + 3\sqrt {69}}{2}\right)^{1/3}.
\end{equation*}
Substituting $m = 2$ and $m = 3$ in (2.13) we get
\begin{align*}
\sum_{k = 1}^\infty \frac{x^{3k}}{k^2\binom{6k}{2k}} &= 12\arctan
^2\left[ {\frac{\sqrt 3 }{2\phi ( - x) - 1}} \right] +
12\arctan ^2\left[ {\frac{\sqrt 3 }{2\phi (x) - 1}}
\right]\\[.4pc]
&\quad\, - \log ^2\left[ {\frac{\phi(- x)^3 + 1}{[\phi ( -
x) + 1]^3}} \right] - \log ^2\left[ {\frac{\phi (x)^3 + 1}{[\phi
(x) + 1]^3}} \right]
\end{align*}
and
\begin{align*}
\sum_{k = 1}^\infty \frac{x^{3k}}{k^2\binom {9k} {3k}} &=
18\arctan ^2\left[ {\frac{\sqrt 3 }{2\phi (x) - 1}} \right] +
18\arctan ^2\left[ {\frac{\sqrt 3 }{2\phi (ax) -
1}}\right]\\[.4pc]
&\quad + 18\arctan ^2\left[ {\frac{\sqrt 3 }{2\phi (a^2x) - 1}}
\right] - \frac{3}{2}\log ^2\left[ {\frac{\phi (x)^3 + 1}{[\phi
(x) + 1]^3}} \right]\\[.4pc]
&\qquad - \frac{3}{2}\log ^2\left[ {\frac{\phi (a x)^3 +
1}{[\phi (a x) + 1]^3}} \right] - \frac{3}{2}\log
^2\left[{\frac{\phi (a^2 x)^3 + 1}{[\phi (a^2 x) + 1]^3}}
\right],
\end{align*}
where $a = (i\sqrt 3 - 1) / 2$.

Of these results, eqs~(3.6), (3.9), (3.11) and (3.12) have been
evaluated by Borwein and Girgensohn \cite{7} experimentally by the
method called {\it integer relation algorithm} which does not
constitute a mathematical proof. So their results are just
conjectural. Our results verify Borwein and Girgensohn's
experimental evaluations. All the results we obtained seem to be
new.

\section*{Acknowledgements}

The author would like to thank the referee for useful suggestions
and for giving the correct version of eqs~(3.1) and (3.2).

\end{document}